\documentclass[12pt]{article}
\pdfoutput=1
\usepackage{amssymb,amsmath,amsthm,latexsym}
\usepackage[cp1251]{inputenc}
\usepackage{dsfont}
\usepackage{amssymb}
\usepackage{indentfirst}
\usepackage{graphicx}
\usepackage{multicol}
\usepackage{amsfonts, array, hhline}
\usepackage{color}
\usepackage{grffile}
\usepackage[export]{adjustbox}
\begin{document}
\title{On a relation between packing and covering densities of convex bodies}
\author{Roman Prosanov \thanks{Supported by the grant from Russian Federation President Programme
of Support for Leading Scientific Schools 6760.2018.1 and SNF grant $200021_-169391$ ``Discrete curvature and rigidity''.}}
\date{}
\AtEndDocument{\bigskip{\footnotesize
\par
  \textsc{Universit\'{e} de Fribourg, Chemin du Mus\'{e}e 23, CH-1700 Fribourg, Switzerland} \par
  \textsc{Moscow Institute Of Physics And Technology, Institutskiy per. 9, 141700, Dolgoprudny, Russia} \par
  \textit{E-mail}: \texttt{rprosanov@mail.ru}
}}
\maketitle

\newtheorem{dfn}{Definition}[section]
\newtheorem{thm}[dfn]{Theorem}
\newtheorem{lm}[dfn]{Lemma}
\newtheorem{crl}[dfn]{Corollary}
\newtheorem{cnj}[dfn]{Conjecture}
\newtheorem{qst}[dfn]{Question}
\renewcommand{\proofname}{Proof}

\renewcommand{\refname}{Bibliography}
\renewcommand{\proofname}{Proof}
\renewcommand{\figurename}{Figure}
\renewcommand{\le}{\leqslant}
\renewcommand{\leq}{\leqslant}
\renewcommand{\ge}{\geqslant}
\renewcommand{\geq}{\geqslant}
\renewcommand{\mathds}{\mathbb}
\newcommand{\e}{\varepsilon}
\newcommand{\R}{\mathbb{R}}
\newcommand{\dist}{{\rm dist}}
\newcommand{\vol}{\textrm{vol}}
\newcommand{\den}{{\rm den}}

\begin{abstract}
We show that a convex body admits a translative dense packing in $\R^d$ if and only if it admits a translative economical covering.
\end{abstract}

\section{Introduction}

\subsection{Packing and covering densities}

Let $C$ be a $d$-dimensional convex body in $\R^d$, i.e. a compact convex set with nonempty interior. A (translative) \emph{arrangement} is a set $C+A$, where $A$ is a discrete point set in $\R^d$. We assume that $A$ is infinite. An arrangement is called \emph{packing} if no two translates of $C$ in $C+A$ have an interior point in common. An arrangement is called \emph{covering} if $\R^d = C+A$.

Define \emph{upper} and \emph{lower} densities of an arrangement
$$\overline{{\rm den}}(C+A) = \limsup_{r \rightarrow \infty} \frac{\sum_{a \in A}\vol\big((C+a)\cap B^d(r)\big)}{\vol(B^d(r))},$$
$$\underline{{\rm den}}(C+A) = \liminf_{r \rightarrow \infty} \frac{\sum_{a \in A}\vol\big((C+a)\cap B^d(r)\big)}{\vol(B^d(r))},$$
where $B^d(r)$ is the Euclidean ball of radius $r$ centered at the origin.

The (translative) \emph{packing density} of $C$ is $$\delta_T(C) = \sup_{C+A {\rm~is~a~packing}} \overline{{\rm den}}(C+A).$$

Similarly, the (translative) \emph{covering density} of $C$ is $$\theta_T(C) = \inf_{C+A {\rm~is~a~covering}} \underline{{\rm den}}(C+A).$$

An important example is a \emph{periodic arrangement}, i.e. an arrangement of the form $C+\Lambda + X$, where $\Lambda$ is a lattice and $X$ is a finite point set. In this case, $$\overline{{\rm den}}(C+\Lambda + X)=\underline{{\rm den}}(C+\Lambda + X) = \frac{|X|\vol(C)}{\vol(\R^d/\Lambda)}.$$ Then we will denote this quantity just as ${\rm den}(C+\Lambda + X)$.

We can consider arrangements consisting not only from translates of $C$, but from any congruent copies of $C$. In this case the packing and covering densities of $C$ can be defined similarly. Denote them by $\delta(C)$ and $\theta(C)$ respectively. Another important case is the case of \emph{lattice arrangements}, i.e. of the form $C+\Lambda$, where $\Lambda$ is a lattice. The corresponding densities over lattice arrangements only are denoted by $\delta_L(C)$ and $\theta_L(C)$.

Bounding packing and covering densities (especially for some specific choices of $C$, e.g. Euclidean balls) is one of the main problems in discrete geometry. Despite a lot of progress, plenty important questions remain open.

Clearly, for any $C$ we have $$\delta_T(C) \leq \delta(C) \leq 1$$ and $$\theta_T(C) \geq \theta(C) \geq 1.$$ The equality of any of these densities to 1 means that copies of $C$ tessellate $\R^d$. Hence, we have $$\delta(C)=1~~~~~~\Longleftrightarrow~~~~~~\theta(C)=1,$$ $$\delta_T(C)=1~~~~~~\Longleftrightarrow~~~~~~\theta_T(C)=1.$$

A natural question arises from this observation: if a body $C$ can not be packed densely, does it mean that it can not cover $\R^d$ economically? In the book by P.~Brass, W.~Moser and J.~Pach~\cite{BMP} this conjecture is attributed to W.~Kuperberg (Conjecture 1 in Chapter 1.10, the original notation is saved):

\begin{cnj}
\label{Kup}
Let $d \geq 2$ be fixed. Then for any $\varepsilon > 0$ there exists a $\delta > 0$ with the property that for every $d$-dimensional convex body $C$, $$(1)~~~~~\delta(C) \leq 1-\varepsilon~~~~~~implies~~~~~~\theta(C) \geq 1 + \delta,$$ $$(2)~~~~~\theta(C) \geq 1+\varepsilon~~~~~~implies~~~~~~\delta(C) \leq 1 - \delta.$$
\end{cnj}

In the notation of Conjecture~\ref{Kup}, $\delta$ and $\delta(C)$ are not the same. The notation $\delta(C)$ is conventional for the packing density and another $\delta$ is common for the epsilon-delta notation. 

The aim of the present note is to prove this conjecture for translative densities.
It will be more helpful to give our statement in the form of converse implications. 

\begin{thm}
\label{main}
Let $d \geq 2$ be fixed.
\vskip+0.2cm
(1a) Let $0< \e \leq \frac{1}{d^{d+1}}$ and $C$ be a $d$-dimensional convex body or ${0<\e<1}$ and $C$ in addition be centrally symmetric. If for the translative packing density we have $\delta_T(C) > 1 - \varepsilon$, then the translative covering density of $C$ satisfies $$\theta_T(C) < \left(1+\e^{\frac{1}{d+1}}\right)^{d+1}.$$
\vskip+0.2cm
(1b) Let $\frac{1}{d^{d+1}} < \e < 1$ and $C$ is not centrally symmetric. If for the translative packing density we have $\delta_T(C) > 1 - \varepsilon$, then the translative covering density of $C$ satisfies $$\theta_T(C) < \left(1+\e d^d\right)\left(1+\frac{1}{d}\right)^d.$$
\vskip+0.2cm
(2) Let $0 < \e <1$ and $C$ be a $d$-dimensional convex body. If for the translative covering density we have $\theta_T(C) < 1 + \varepsilon$, then the translative packing density of $C$ satisfies $$\delta_T(C) > \left(1-\e^{\frac{1}{d+1}}\right)^{d+1}.$$
\end{thm}

Clearly, Conjecture~\ref{Kup} for translative densities follows from this theorem.


\subsection{Previous results and future perspectives}

The only already known case of Conjecture~\ref{Kup} was established by Ismailescu in~\cite{Ism}. He considered $d=2$ and only centrally symmetric convex bodies. More precisely, he showed that in this case $$1-\delta_L(C) \leq \theta_L(C) -1 \leq 1.25\sqrt{1-\delta_L(C)}.$$

L.~Fejes T\'oth~\cite{FeT} established that for any planar centrally symmetric convex body $C$ we have $\delta(C)=\delta_T(C)=\delta_L(C)$ and Dowker~\cite{Dow} proved that $\theta_L(C) = \theta_T(C)$. 
Hence, Ismailescu's result extends to the case of all translative arrangements. His proof is based on the approximation of $C$ by centrally symmetric octagons and can not be extended to higher dimensions.

In~\cite{FTK} G. Fejes T\'oth and W.~Kuperberg proposed to understand links between packing and covering densities in a more general way. They defined the set $\Omega_d$ (resp. $\Omega^*_d$) of points $(x, y) \in \R^2$ such that there exists a $d$-dimensional convex (resp. centrally symmetric) body $C$ with $\delta(C)=x$ and $\theta(C)=y$. The definition of $\Omega_d$ (resp. $\Omega^*_d$) can be restricted to the case of translative or lattice densities. It may be of interest to characterize these sets. In fact, it is still unknown whether these sets are closed (but it is known for translative or lattice cases) or convex. We refer the reader to the paper~\cite{Ku13} investigating the planar case. Several inequalities involving both $\delta_L(C)$ and $\theta_L(C)$ were established e.g. in~\cite{IK},~\cite{Ku87} and~\cite{Sm}, but also only in low dimensions.

In order to prove the translative Kuperberg conjecture it is enough to consider only sufficiently small values of $\e$ with respect to $d$. When $\e$ is not very small and $d$ is sufficiently large it is interesting to compare Theorem~\ref{main} with the best known general bounds on packing and covering densities. 


In the case of coverings by translates of an arbitrary $d$-dimensional convex body $C$ the following inequality was established by G.~Fejes T\'oth~\cite{FT} (which slightly improves a previous result by Rogers): $$\theta_T(C) \leq d\ln d+ d\ln\ln d +o(d).$$

We see that Theorem~\ref{main} gives us a stronger bound if  $$1-\delta_T(C) <\frac{\ln d}{ed^{d-1}}$$ or if $$1-\delta_T(C)  < \left(\frac{\ln (d \ln d +d\ln\ln d)}{d+1}\right)^{d+1}$$ and $C$ in addition be centrally symmetric.

For packing densities of centrally symmetric convex bodies and $d$ sufficiently large the following result by Schmidt~\cite{Sch} is the best known: $$\delta_T(C)\geq \delta_L(C) \geq\frac{cd}{2^d}.$$

Comparing with the last inequality in Theorem~\ref{main} we see that the latter gives a stronger bound if $$\theta_T(C)-1< \left(\frac{1}{2}- \frac{\ln(2cd)}{d+1}\right)^{d+1}.$$

For a non centrally symmetric $C$ we should use the observation of Minkowski (see \cite{HB}, Chapter~2): $$\delta_T(C)=2^d\delta_T(C-C)\frac{\vol(C)}{\vol(C-C)}.$$

Denote $\left(\frac{\vol(C-C)}{\vol(C)}\right)^{1/d}$ by $W(C)$. Then we have $$\delta_T(C) \geq \frac{cd}{W^d(C)}.$$


Hence, we obtain a better bound provided $d$ sufficiently large and $$\theta_T(C)-1< \left( 1 - \frac{1}{W(C)}-\frac{\ln\left(W(C)cd\right)}{d+1}\right)^{d+1}.$$

An interesting question is to understand if the dependence of our bounds from $d$ is necessary. In other words, we would like to propose the following problem:

\begin{qst}
Is it true that for any $\varepsilon > 0$ there exists $\mu > 0$ with the property that for every $d$ and every $d$-dimensional convex body $C$, $$(1)~~~~~\delta(C) \leq 1-\varepsilon~~~~~~implies~~~~~~\theta(C) \geq 1 + \mu,$$ $$(2)~~~~~\theta(C) \geq 1+\varepsilon~~~~~~implies~~~~~~\delta(C) \leq 1 - \mu .$$
\end{qst}

It is also of interest to prove an analogue of Theorem~\ref{main} for the case of lattice densities only (it is conjectured~\cite{BMP} that in higher dimensions $\delta_L(C) \neq \delta_T(C)$ and $\theta_L(C) \neq \theta_T(C)$). Such a proof should use totally different ingredients.

Our proof is based on a bound on the number of steps of a certain greedy algorithm. It seems that in previous years most of results on packing and covering densities in higher dimensions used averaging or probabilistic arguments. Recently the focus started to shift in the direction of more deterministic techniques. For instance, in~\cite{Na} Nasz\'odi gave a new proof of some well-known covering results via discretization and a lemma connecting fractional covering numbers of finite hypergraphs with integral ones. Proofs of this lemma considered a greedy algorithm applied to finite sets. In~\cite{RV} Rolfes and Vallentin explored a greedy algorithm applied directly to geometric covering problems and also obtained several classical results through their method. 

Packing and covering results have some applications to other problems in discrete geometry. For instance, consider the problem of finding \emph{the chromatic number $\chi(S)$} of a subset $S \subseteq \R^d$. The chromatic number $\chi(S)$ is the minimal number of colors sufficient to color $S$ in such a way that any two points at the distance 1 have different colors. Using deterministic covering algorithms in~\cite{Pr1} the author gave a new proof of the upper bound for $\chi(\R^d)$ and in~\cite{Pr2} the author established new upper bounds for $\chi(S^{d-1}_R)$, where $S^{d-1}_R$ is a $(d-1)$-dimensional Euclidean sphere of radius $R$. For more details about geometric chromatic numbers the reader is refereed to the surveys~\cite{Ra2},~\cite{Ra1}.

\textbf{Acknowledgments.} The author would like to thank I. Izmestiev and A. Polyanskii for useful discussions and remarks.

\section{Proof of Theorem~\ref{main}}










We need an important theorem of Rogers (see~\cite{Ro}, Theorems 1.7 and 1.9).

\begin{thm}
\label{Ro}
For a convex body $C$ in the definition of its translative packing density we can take the supremum over only periodic arrangements. The same holds for its translative covering density.
\end{thm}

We assume that $C$ contains the origin in the interior. By $\lambda C + x$ we denote the image of $C$ under the composition of the homothety with the center at the origin and the coefficient $\lambda$ and the translation by the vector $x$.

Now we are able to give a proof of the main theorem.

\textit{Proof of (1).} Assume that $\delta_T(C) > 1 - \e$. By Theorem~\ref{Ro} there is a lattice $\Lambda$ and a finite point set $X$ such that $C+\Lambda+X$ is a packing and $$\den(C+\Lambda+X) > 1 - \e.$$ Consider the torus $T=\R^n/\Lambda$. The sets $X$ and $C$ can be projected to $T$. Abusing the notation, we still denote by $X$ and $C$ their images under this projection. This will not lead to an ambiguity as from now on we work only on $T$. The arrangement $C+X$ is a packing on $T$. 

Let $k=|X|$. Then $$\frac{k\vol(C)}{\vol(T)} = \den(C+\Lambda+X) > 1- \e.$$ As our problem is homothety invariant, we may assume that $\vol(T)=1$. Hence, $$k\vol(C) > 1-\e.$$

Define $S_0=\vol(T \backslash (C+X))$. We have $S_0 < \e$. Also, let $X_0=X$.

Fix $0< \alpha < 1$. We proceed iteratively. Assume that $\left(1+\alpha \right)C + X_i$ does not cover $T$. Then there exists $y \in T$ which is not covered by $\left(1+\alpha\right)C+X_i$. We have that for every $x \in X$, $$\left(-\alpha C + y\right) \cap \left(C + x\right) = \emptyset.$$
Indeed, $-\alpha C + y$ can be obtained as the image of $C + x$ under a homothety with the coefficient $-\alpha$ and the center in the segment $xy$ laying outside of $C + x$. Next, we are looking for $y'$ such that $C+y'$ covers $-\alpha C + y$. If $C$ is centrally symmetric, then we can take $y'=y$. In the other case we need the condition $\alpha \leq \frac{1}{d}$. Then the existence of such $y'$ follows from the fact that we can put a translate of $-\frac{1}{d}C$ into $C$. This statement is equivalent to the existence of a point in the interior of $C$ such that every chord through this point is divided in the ratio not greater than $d$. It is a well-known implication of the Helly theorem and was proved by Minkowski and Radon, see e.g.~\cite{To}, Corollary 1.4.2.

Consider $X_{i+1} = X_i \cup \{y'\}$ and $S_{i+1} = \vol(T \backslash (C+X_{i+1}))$. Clearly, $$0 \leq S_{i+1} \leq S_i - \vol\left(-\alpha C\right) = S_i - \alpha^d \vol(C) < \e - (i+1)\alpha^d\vol(C).$$  If $\left(1+\alpha\right)C + X_{i+1}$ does not cover $T$, then repeat this process. At every step $S_i \geq 0$. Hence, we stop after $l$ steps and the upper bound on $l$ can be deduced from the inequality $$\e > l \alpha^d \vol(C).$$ We rewrite it as $$l\vol(C) < \frac{\e}{\alpha^d}.$$

We obtain that $\left(1+\alpha\right)C + X_l$ covers $T$. Then $\left(1+\alpha\right)C + X_l + \Lambda$ covers $\R^d$. Now we need to estimate the density of this arrangement $$\theta_T(C)=\theta_T(\left(1+\alpha\right)C)\leq$$ $$ \leq \den\left(\left(1+\alpha \right)C+X_l+\Lambda\right) = |X_l| \vol\left(\left(1+\alpha \right)C\right) =$$ $$= (k+l)\left(1+\alpha\right)^d\vol(C) = (k\vol(C)+l\vol(C))\left(1+\alpha\right)^d.$$

As $C+X_0$ is a packing, then $k\vol(C) \leq 1$. Using this and $l\vol(C) < \frac{\e}{ \alpha^d}$ we get $$\theta_T(C) < \left(1+\frac{\e}{ \alpha^d}\right)(1+\alpha)^d.$$

After the calculation of the derivative in $\alpha$ we can see that for fixed $d$ and $\e$ this expression attains its minimal value at $\alpha=\e^{\frac{1}{d+1}}$. If $\e^{\frac{1}{d+1}} \leq \frac{1}{d}$, then it is an admissible value for any $C$. After the substitution we obtain the bound $$\theta_T(C) < \left(1+\e^{\frac{1}{d+1}}\right)^{d+1}.$$

If $\e^{\frac{1}{d+1}} > \frac{1}{d}$ and $C$ is not centrally symmetric, then the admissible value of $\alpha$ minimizing the expression at the right-hand side is $\alpha=\frac{1}{d}$. In this case we have $$\theta_T(C) < \left(1+\e d^d\right)\left(1+\frac{1}{d}\right)^d.$$



\qed

\textit{Proof of (2).} Assume that $\theta_T(C) < 1 + \e$. Similarly, by Theorem~\ref{Ro} there is a lattice $\Lambda$ and a finite point set $X$ such that $C+\Lambda+X$ is a covering and $$\den(C+ \Lambda+X) < 1 + \e.$$

Moreover, we can choose $\Lambda$ such that for any $\lambda_1$ and $\lambda_2 \in \Lambda$, the translate $C+\lambda_1$ does not intersect $C+\lambda_2$. Indeed, let $m > 0$ be an integer. For every element $\gamma$ of the group $\Lambda / m\Lambda$ choose a representative $\lambda(\gamma) \in \Lambda$. By $\Lambda_m \subset \R^d$ denote the set $\{\lambda(\gamma) : \gamma \in \Lambda/m\Lambda\}$. There exists $m$ such that the desired condition is satisfied for $m\Lambda$. Take $X' = X + \Lambda_m$. Then $C+m\Lambda+X'$ is a covering and $$\den(C+m\Lambda+X') = \den(C+\Lambda+X) < 1+ \e.$$ Then we can replace $\Lambda$ with $m\Lambda$ and $X$ with $X'$.

Let $T$ be the torus $\R^n/\Lambda$ and $k=|X|$. As in the proof of (1) from now on we consider $X$ and $C$ as subsets of $T$. Then $C+X$ is a covering of $T$ and $$\frac{k\vol(C)}{\vol(T)} = \den(C+\Lambda+X) < 1+ \e.$$ As previously, assume that $\vol(T)=1$. Hence, $$k\vol(C) < 1+\e.$$

Define $S_0=k\vol(C)-1$. Then $S_0 < \e$. Also, let $X_0=X$.

Fix $0< \alpha < 1$. Now we proceed iteratively. Assume that $\left(1-\alpha \right)C + X_i$ is not a packing. Then there exists $y \in T$ and $x, x' \in X_i$ such that $$\left(\alpha C + y\right) \subset (C+x)\cap(C+x').$$ Indeed, there are $x, x' \in X_i$ such that $$\left(\left(1-\alpha\right)C+x\right)\cap\left(\left(1-\alpha\right)C+x'\right) \neq \emptyset.$$ Choose $y$ in their intersection. Then clearly $$\left(\alpha C + y\right) \subset \left(\alpha C + \left(1-\alpha\right)C+x \right) = C+x.$$ Similarly, $\left(\alpha C + y\right) \subset  C+x'$.

Consider $X_{i+1} = X_i \backslash \{x'\}$ and $$S_{i+1} =|X_{i+1}|\vol(C)-(1-\vol(T \backslash (C+X_{i+1}))).$$ Naturally, $S_i$ measures the covering excess of the arrangement $C+X_i$. We obtain $$0 \leq S_{i+1} \leq S_i - \vol\left(\alpha C\right) = S_i - \alpha^d\vol(C) < \e - (i+1)\alpha^d \vol(C).$$  If $\left(1-\alpha \right)C + X_{i+1}$ is not a packing, then repeat this process. At every step $S_i \geq 0$. Hence, as previously we will stop after $l$ steps, where $$l\vol(C) < \frac{\e}{\alpha^d}.$$

We have $\left(1-\alpha\right)C + X_l$ is a packing in $T$ and $\left(1-\alpha\right)C + X_l+\Lambda$ is a packing in $\R^d$. Now we need to estimate the density of this arrangement $$\delta_T(C)=\delta_T(\left(1-\alpha\right)C)\geq \den\left(\left(1-\alpha\right)C + X_l +\Lambda\right) =$$ $$= |X_l| \vol\left(\left(1-\alpha\right)C\right) = (k-l)\left(1-\alpha\right)^d\vol(C) =$$ $$= (k\vol(C)-l\vol(C))\left(1-\alpha\right)^d >  \left(1-\frac{\e}{\alpha^d}\right)\left(1-\alpha\right)^d.$$

This expression is minimized as $\alpha = \e^{\frac{1}{d+1}}$. Then we obtain $$\delta_T(C)> \left(1-\e^{\frac{1}{d+1}}\right)^{d+1}.$$


\qed

\bibliographystyle{abbrv}
\bibliography{Kuperberg}

\begin{thebibliography}{10}

\bibitem{BMP}
P.~Brass, W.~Moser, and J.~Pach.
\newblock {\em Research problems in discrete geometry}.
\newblock Springer, New York, 2005.

\bibitem{Dow}
C.~H. Dowker.
\newblock On minimum circumscribed polygons.
\newblock {\em Bull. Amer. Math. Soc.}, 50:120--122, 1944.

\bibitem{FT}
G.~Fejes~T\'{o}th.
\newblock A note on covering by convex bodies.
\newblock {\em Canad. Math. Bull.}, 52(3):361--365, 2009.

\bibitem{FTK}
G.~Fejes~T\'oth and W.~Kuperberg.
\newblock A survey of recent results in the theory of packing and covering.
\newblock In {\em New trends in discrete and computational geometry}, volume~10
  of {\em Algorithms Combin.}, pages 251--279. Springer, Berlin, 1993.

\bibitem{FeT}
L.~Fejes~T\'{o}th.
\newblock Some packing and covering theorems.
\newblock {\em Acta Sci. Math. Szeged}, 12:62--67, 1950.

\bibitem{HB}
J.~E. Goodman, J.~O'Rourke, and C.~D. T\'{o}th, editors.
\newblock {\em Handbook of discrete and computational geometry}.
\newblock Discrete Mathematics and its Applications. CRC Press, Boca Raton, FL,
  2018.

\bibitem{Ism}
D.~Ismailescu.
\newblock Inequalities between lattice packing and covering densities of
  centrally symmetric plane convex bodies.
\newblock {\em Discrete Comput. Geom.}, 25(3):365--388, 2001.

\bibitem{IK}
D.~Ismailescu and B.~Kim.
\newblock Packing and covering with centrally symmetric convex disks.
\newblock {\em Discrete Comput. Geom.}, 51(2):495--508, 2014.

\bibitem{Ku87}
W.~Kuperberg.
\newblock An inequality linking packing and covering densities of plane convex
  bodies.
\newblock {\em Geom. Dedicata}, 23(1):59--66, 1987.

\bibitem{Ku13}
W.~Kuperberg.
\newblock The set of packing and covering densities of convex disks.
\newblock {\em Discrete Comput. Geom.}, 50(4):1072--1084, 2013.

\bibitem{Na}
M.~Nasz\'odi.
\newblock On some covering problems in geometry.
\newblock {\em Proc. Amer. Math. Soc.}, 144(8):3555--3562, 2016.

\bibitem{Pr1}
R.~{Prosanov}.
\newblock {A new proof of the Larman-Rogers upper bound for the chromatic
  number of the Euclidean space}.
\newblock {\em ArXiv e-prints}, Oct. 2016.

\bibitem{Pr2}
R.~Prosanov.
\newblock Chromatic numbers of spheres.
\newblock {\em Discrete Math.}, 341(11):3123--3133, 2018.

\bibitem{Ra2}
A.~{Raigorodskii}.
\newblock {Combinatorial geometry and coding theory.}
\newblock {\em {Fundam. Inform.}}, 145(3):359--369, 2016.

\bibitem{Ra1}
A.~M. {Raigorodskii}.
\newblock {Coloring distance graphs and graphs of diameters.}
\newblock In {\em {Thirty essays on geometric graph theory}}, pages 429--460.
  Berlin: Springer, 2013.

\bibitem{Ro}
C.~A. Rogers.
\newblock {\em Packing and covering}.
\newblock Cambridge Tracts in Mathematics and Mathematical Physics, No. 54.
  Cambridge University Press, New York, 1964.

\bibitem{RV}
J.~H. Rolfes and F.~Vallentin.
\newblock Covering compact metric spaces greedily.
\newblock {\em Acta Math. Hungar.}, 155(1):130--140, 2018.

\bibitem{Sch}
W.~M. Schmidt.
\newblock On the {M}inkowski-{H}lawka theorem.
\newblock {\em Illinois J. Math.}, 7:18--23, 1963.

\bibitem{Sm}
E.~Smith.
\newblock An improvement of an inequality linking packing and covering
  densities in 3-space.
\newblock {\em Geom. Dedicata}, 117:11--18, 2006.

\bibitem{To}
G.~Toth.
\newblock {\em Measures of symmetry for convex sets and stability}.
\newblock Universitext. Springer, Cham, 2015.

\end{thebibliography}

\end{document}